\documentclass[letterpaper, 10 pt, conference]{ieeeconf}  

\IEEEoverridecommandlockouts                              

\usepackage{cite}
\usepackage{amsmath,amssymb,amsfonts}
\usepackage{algorithmic}
\usepackage{graphicx}
\usepackage{textcomp}
\usepackage{xcolor}
\usepackage{svg}
\usepackage{graphicx}
\usepackage{booktabs}
\usepackage{multirow}
\usepackage{epstopdf}          
\usepackage{stfloats}
\usepackage{subcaption}
\usepackage{hyperref}
\usepackage{xcolor}

\definecolor{linkblue}{RGB}{0, 150, 255} 

\hypersetup{
    colorlinks=true,
    linkcolor=blue,
    filecolor=magenta,      
    urlcolor=linkblue,
}
\DeclareGraphicsExtensions{.pdf,.png,.jpg,.eps}

\title{\LARGE \bf
Reference-Free, Long-Horizon Trajectory Optimization \\for Aggressive Autonomous Driving in Milliseconds
}

\author{Prayag Sharma$^{1}$, Jonathan Y.M. Goh$^{2}$, and Franck Djeumou$^{1}$
\thanks{*This work was funded by Toyota Research Institute (TRI)}
\thanks{$^{1}$P. Sharma and F. Djeumou are with MANE Dept., Rensselaer Polytechnic Institute, Troy, NY, USA
        {\tt\small {sharmp6 and djeumf2}@rpi.edu}}%
\thanks{$^{2}$ J. Goh is with Toyota Research Institute, Los Altos, CA, USA
        {\tt\small jon.goh@tri.global}}%
}

\begin{document}

\maketitle
\thispagestyle{empty}
\pagestyle{empty}

\begin{abstract}
Autonomous vehicles must generate long-horizon and dynamically feasible trajectories in real time—even when operating at the limits of vehicle handling—to ensure safe operation in adverse conditions. However, existing work rarely quantifies the computational demands of generating such trajectories without prior references, warm starts and often defaults to low-fidelity models, compromising accuracy and control authority. We investigate the modeling and solver design choices that enable real-time solution of long-horizon, reference-free optimal control problems (OCPs) using full vehicle dynamics.
To this end, we analyze vehicle stiffness properties to justify the OCP's integration scheme and show that lower-order A-stable methods consistently outperform alternatives, with solve time differences reaching two orders of magnitude. We show that robust nonlinear solver performance hinges on understanding barrier parameter update strategies and safeguarding techniques for Hessian indefiniteness, inherent in some interior point methods. Lastly, we propose a computationally efficient method for generating initial guesses using dynamic equilibrium, unlocking real-time performance and reducing initial infeasibility by up to four orders of magnitude. Extensive benchmarking and high-fidelity \textrm{BeamNG} simulation demonstrate compute times as low as 55 ms over a 260 m horizon, including high-speed obstacle avoidance scenarios where drifting emerges as a necessary component of feasible trajectory generation.

\end{abstract}

\section{INTRODUCTION}


To ensure safety in all critical situations, autonomous vehicles must be engineered to exploit the absolute limits of their dynamic capabilities, pushing beyond the boundaries of normal vehicle operation. Existing work has established that intentionally pushing a vehicle beyond traditional operational limits not only expands the safety envelope but can often be the only feasible evasive maneuver \cite{CA_scenarios, justifying_emergency, Lund}. Such maneuvers require long-horizon and dynamically feasible planning for linking immediate actions to their downstream consequences in unforeseen situations. Unlocking this capability requires a paradigm shift from conventional short-horizon reference tracking to solving a \textit{reference-free, long-horizon} optimal control problem (OCP) in real time (5-10 Hz). This paper presents a systematic investigation into the fundamental design choices required to construct such a framework, thereby expanding the operational envelope of autonomous vehicles.


Prior approaches underscore the complexity of solving the full nonlinear OCPs from scratch, identifying the underlying problems as highly sensitive to initial guesses and too computationally intensive for real-time deployment\cite{CA_scenarios, Lund, switch_logic, motion_primitives, long_horizon_point_mass}. This has led to a shared strategy across both high-performance racing \cite{subosits2019racetrack,thompson2024adaptive, chris_g_cascaded_22} and autonomous drifting \cite{goh2019_thesis,goh2024beyond,weber2023modeling}: decoupling the problem into offline reference generation and online tracking. Even though they are successful in handling even pop-up obstacles \cite{subosits2019racetrack,thompson2024adaptive}, the core dependency of such approaches to an offline-generated reference precludes their adaptation to a dynamically changing course. Recent attempts toward real-time generation include a task-specific OCP to transition between drift equilibria \cite{talbot2024optimal}, and reinforcement learning (RL) to generate drift trajectories \cite{djeumou2025reference, domberg2022deep}. However, the former remains highly task-specific, while RL's reliance on training data makes its reliability in unfamiliar scenarios difficult to guarantee.

\begin{figure}[!t]
  \centering
  \includegraphics[width=0.9\columnwidth]{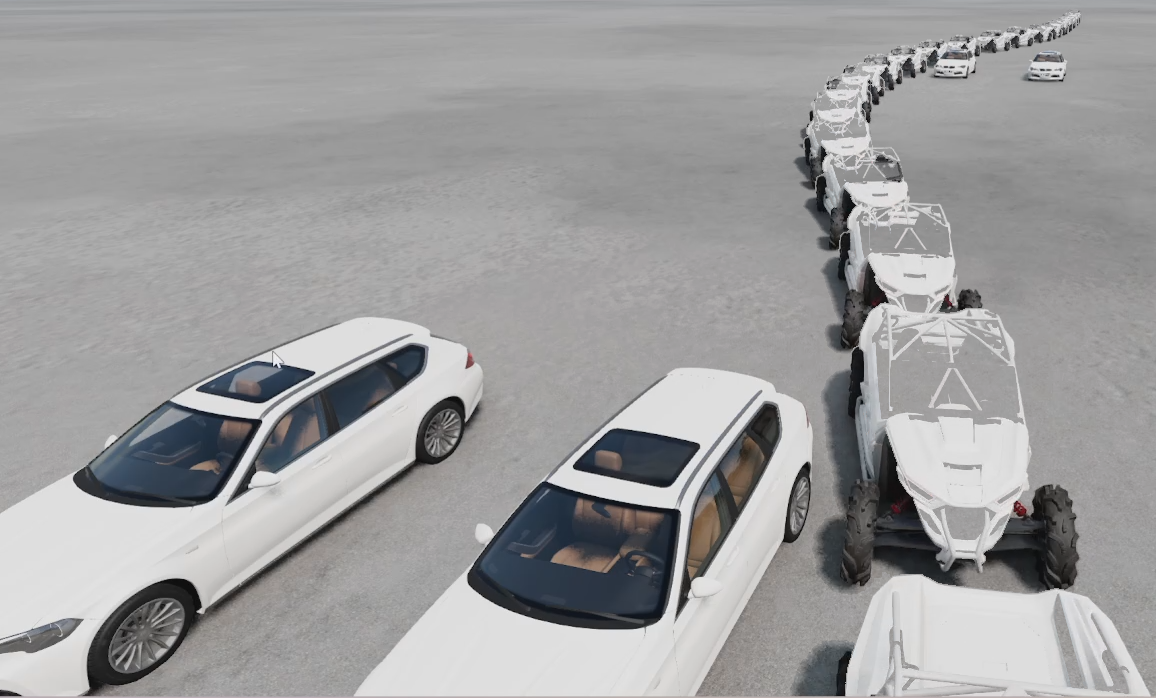}
  \caption{Tracking of high-speed emergency collision avoidance trajectories in BeamNG simulation. Video available at  \protect\href{https://drive.google.com/file/d/1QN7d9_mewr0E-bx7v-Bl92dXdkhoglsv/view?usp=sharing}{\color{linkblue}https://tinyurl.com/trajopt}}
  \label{fig:beamNG}
  \vspace{-6mm}
\end{figure}



To enable online replanning in safety-critical scenarios, many existing approaches reduce model fidelity by solving online long-horizon Optimal Control Problems using simplified point-mass dynamics \cite{long_horizon_point_mass}. Although computationally tractable, these abstractions often produce trajectories that violate dynamic feasibility \cite{motion_primitives}. Alternatively, systems based on pre-computed motion primitives \cite{motion_primitives,switch_logic} constrain the vehicle to a fixed library of maneuvers, limiting expressiveness and proving brittle in scenarios not anticipated during design. Both workarounds trade off solving the underlying OCPs for solutions that may either result in dynamic infeasibility or may fail to exploit the full vehicle capabilities.

Although existing work excels at reference tracking, the foundational challenge of real-time reference generation from scratch remains largely unaddressed, favoring offline or task-specific solutions. As a consequence, critical design choices that directly impact real-time performance, such as an accurate and efficient numerical integrator, the nonlinear solver strategy, and a strategy for a generalizable initial guess, are rarely examined. This paper addresses this gap by presenting a systematic design of a trajectory optimization framework for solving the full nonlinear, reference-free OCP in real time. Our contributions are as follows:

\noindent\textbullet~Through extensive benchmarking, we show that, while interior point methods \cite{wachter2006ipopt,byrd2006knitro} outperform sequential quadratic programming (SQP) approaches\cite{gill2005snopt} for trajectory generation, commonly used solvers like IPOPT \cite{wachter2006ipopt} are unsuitable for real-time use. Our findings reveal that robust and real-time performance emerges from carefully selected barrier update and numerical ill-conditioning handling strategies found in solvers like KNITRO\cite{byrd2006knitro}.\\
\noindent\textbullet~We analyze the stiffness properties of vehicle dynamics and establish a principled basis for selecting stable and computationally efficient integration schemes. Further, our evaluation of nine numerical integrators reveals computational performance differences of up to two orders of magnitude. Critically, we demonstrate that due to system stiffness, lower-order A-stable methods surprisingly outperform higher-order, non-A-stable methods in terms of accuracy.\\
\noindent\textbullet~ We propose a strategy to generate a high-quality, reusable initial guess for the OCP by solving a single-node, box-constrained least-square problem formulation. Such a guess is computationally cheap to obtain, is by construction dynamically feasible, reduces total initial feasibility error for the optimizer by three to four orders of magnitude compared to an ill-informed guess, while enabling real-time performance. We then show the robustness of this guess by solving OCPs spanning both time-optimal racing and collision avoidance.\\
\noindent\textbullet~ We validate our approach on a critical collision avoidance scenario, where drifting around the obstacles is the only feasible solution. Our framework generates a drifting trajectory at highway speeds ($90$ km/hr) in real time ($72$ ms), while constraining emerging drift behavior renders the problem infeasible.
To confirm real-world viability, we demonstrate successful tracking of trajectories for both racing and obstacle avoidance in a high-fidelity BeamNG\cite{beamng_tech} simulation.

\section{Problem Formulation}
\subsection{Vehicle Dynamics}
We use a single-track vehicle model in curvilinear coordinates\cite{goh2019_thesis}, considering the wheel speed and load transfer dynamics in a curvilinear coordinate system\cite{weber2023modeling}. The vehicle's position is described by its arc length, $s$, and lateral error, $e$, relative to a reference path with curvature $k_{\mathrm{ref}}(s)$. The equations of motion are given by:
\begin{equation}
\label{eq:dynamics_full}
\frac{d}{ds}
\begin{bmatrix}
    r \\
    V \\
    \beta \\
    V_{\omega r} \\
    \Delta F_z \\
    e \\
    \Delta\psi \\
    s\\
    t
\end{bmatrix}
=
\frac{1}{\dot{s}}
\begin{bmatrix}
    (a (F_{xf} \sin\delta + F_{yf} \cos\delta) - b F_{yr})/I_z\\ 
    (F_{xf} \cos(\delta-\beta) - F_{yf} \sin(\delta-\beta) + \\F_{xr} \cos\beta + F_{yr} \sin\beta)/m\\
    -r + (F_{xf} \sin(\delta-\beta) + F_{yf} \cos(\delta-\beta)\\ - F_{xr} \sin\beta + F_{yr} \cos\beta)/(mV) \\
    R_w(T_{\mathrm{comb}} - R_w F_{xr})/I_w \\
    -c_L (\Delta F_z - h_{cg} F_{x,net}/L) \\
    V \sin(\Delta\psi) \\
    \dot\beta + r - k_{\mathrm{ref}}(s) \dot{s} \\
    1\\
    1-k_{\mathrm{ref}}(s)e/(V\cos(\Delta\psi))
\end{bmatrix}
\end{equation}
where the state $\mathbf{x} = [{r}, {V}, \beta, {V}_{\omega r}, \Delta{F}_z, {e}, \Delta\psi, {s}, t]^T$, the control $\mathbf{u} = [\delta,T_{\mathrm{comb}}, T_{\mathrm{bf}}]^T$, and the independent variable $s$ is such that $\dot{s} = V\cos(\Delta\psi)/(1-k_{\mathrm{ref}}(s)e)$.
Here, $\Delta\psi$ is the velocity vector orientation error w.r.t. path, $r$, $V$, $\beta$ are the yaw rate, vehicle speed, and side slip angle respectively, $V_{\omega r}$ is the rear wheel longitudinal speed, $\Delta F_z$ is the dynamic load transfer (front to rear), $F_{x,net} = F_{xr} + F_{xf}\cos\delta - F_{yf}\sin\delta$ is the net longitudinal force, and $\delta$ is steering angle. We combine the engine torque $T_\mathrm{\mathrm{eng}}\ge0$ and rear brake torque $T_{\mathrm{br}}\le0$ to $T_{\mathrm{comb}}$ such that $T_{\mathrm{eng}}=\max(0,T_{\mathrm{comb}})$ and $T_{\mathrm{br}}=\min(0,T_{\mathrm{comb}})$. As for the model parameters, $m$ is the vehicle mass, $I_z$ the inertia, $a, b, L, h_{cg}$ geometric parameters of the vehicle (front and rear distances from the center of mass, wheelbase, CG height), and $R_w, I_w$ are the wheel radius and rear wheel inertia, respectively.

To model tire forces $(F_{xf},F_{yf},F_{xr},F_{yr})$, we use the isotropic coupled slip brush Fiala model \cite{goh2019_thesis,weber2023modeling}
\begin{align*}
F_{yf} &= -F_{\textrm{total},f}\tan(\alpha_f)/\sigma_f, 
&&F_{xf} = r_\textrm{w}T_{\textrm{br}}, 
\\
F_{yr} &= -F_{\textrm{total},r}\tan(\alpha_r)/\sigma_r,
&&F_{xr} = F_{\textrm{total},r}\kappa_r/\sigma_r.
\label{eq:forces}
\end{align*}
The magnitude of the tire forces $(F_{\textrm{total},f},F_{\textrm{total},r})$ are
\begin{align*}
F_\textrm{total} &=  
\begin{cases}
C\sigma - \frac{C^2 \sigma^2}{3F_\textrm{max}}  + \frac{C^3\sigma^3}{27(F_\textrm{max})^2} &\text{if } |\sigma| < \sigma_\textrm{sl}\\
F_\textrm{max} &\text{if } |\sigma| \geq \sigma_\textrm{sl}\\
\end{cases},
\end{align*}
where 
$(\sigma_f,\sigma_r)$ are the total tire slips, and $(\sigma_{\textrm{sl},f},\sigma_{\textrm{sl},r})$ are the total slips as the tires begin fully sliding
\begin{align*}
\sigma = \sqrt{\tan(\alpha)^2 + \kappa^2},
\quad
\sigma_\textrm{slip}=\tan^{-1}(3\mu F_z/ C).
\end{align*}
The tire loads $(F_{zf},F_{zr})$ depend on the static tire loads $(F_{\textrm{nom},zf}, F_{\textrm{nom},zr})$ as $F_{zf}=F_{\textrm{nom},zf}-\Delta F_z$ and $F_{zr}=F_{\textrm{nom},zr}+\Delta F_z$. 
The maximal tire forces $F_{\textrm{max}}$ are
$$
F_{\textrm{max},f}=\sqrt{(\mu F_{zf})^2-(r_\textrm{w}T_{{\textrm{br}}})^2}.
\quad
F_{\textrm{max},r}=\mu F_{zr},
$$
The slip angles $(\alpha_f,\alpha_r)$ and slip ratios $(\kappa_f,\kappa_r)$ are 
{\small
\begin{align*}
\tan(\alpha_f+&\delta) = (V\sin\beta+ar)/(V\cos\beta),
\ \ 
\\[1mm]
\tan(\alpha_r)&=\frac{V\sin\beta-br}{V\cos\beta}, \ \ 
\kappa_{r,f}= \frac{V_{\omega_{r,f}}-V\cos\beta}{V\cos\beta}.
\end{align*}
}%

\subsection{Minimum Time Discrete Optimal Control problem}
We employ direct multiple shooting with both states and control as the decision variables of our optimization problem. The state and control variables are discretized for a constant $\Delta s_k := s_{k+1}-s_k$. $\mathcal{F}$ is defined as a general map $\mathcal{F}: \mathbb{R}^{n_x} \times \mathbb{R}^{n_x} \times\mathbb{R}^{n_u} \times \mathbb{R} \to \mathbb{R}^{n_x}$ that approximates the next state $x_{k+1}$. Defining the lower and upper bounds on control $\underline{u},\overline{u}$ and its rate $\underline{\dot{u}},\overline{\dot{u}}$, the OCP is formulated as (\ref{eq:dynamics_full}).
We define a common objective of minimum time for both racing and collision avoidance studies as it serves well for both. Except for the variables fixed at boundary conditions(e.g., $s$, $t$, $e$), the above problem is a free initial $\mathbf{x}_{0}$ and final state $\mathbf{x}_{f}$ problem, which helps maintain a general structure. In addition to the track boundaries, we impose bounds on the vehicle's velocity, as this was empirically found to accelerate solver convergence. Furthermore, vehicle-specific control bounds and slew rate constraints are also enforced. The tire saturation constraints prevent the wheels from going beyond their slip limits for racing and can be relaxed for collision avoidance scenarios to allow full dynamic range. A small constant $\epsilon$ is introduced in the tire saturation constraints to prevent infeasibility from brief, instantaneous saturation.
\begin{equation}\label{eq:ocp}
\begin{aligned}
&\min_{\{{x}_k\}_{k=0}^{N},\{u_k\}_{k=0}^{N-1}} \quad t_N \\[-1pt]
&\text{s.t.}\quad
x_{k+1}=\mathcal{F}(x_{k+1},x_k,u_k;\Delta s_k), ~~ k=0,\ldots,N-1, \\
& \text{\scriptsize(initial)}\ \left\{
\begin{aligned}
& s_0=0,\ t_0=0,\ \Delta\psi_0=0,\\& e_{\min}\le e_0\le e_{\max}
\end{aligned}\right. \\
& \text{\scriptsize(track/path bounds)}\ \left\{
\begin{aligned}
& e_{\min}\le e_k\le e_{\max}, \\&V_{\min}\le V_k\le V_{\max}
\end{aligned}\right. ~~ k=0,\ldots,N \\
& \text{\scriptsize(control bounds)}\ \left\{
\begin{aligned}
& \underline{u}\ \le\ u_k\ \le\ \overline{u}
\end{aligned}\right. ~~ k=0,\ldots,N-1 \\
& \text{\scriptsize(control rates)}\ \left\{
\begin{aligned}
& \underline{\dot{u}}\ \le\ \dfrac{u_k-u_{k-1}}{\,t_k-t_{k-1}\,}\ \le\ \overline{\dot{u}}
\end{aligned}\right. ~~ k=1,\ldots,N-1 \\
& \text{\scriptsize(tire saturation)}\ \left\{
\begin{aligned}
&\sigma_{f,sl} - |{\sigma_{f}}_{k}| < \epsilon,\\&\sigma_{r,sl} - |{\sigma_{r}}_{k}| < \epsilon
\end{aligned}\right. ~~ k=0,\ldots,N-1 \\
& \text{\scriptsize(terminal)}\ \left\{
\begin{aligned}
& s_N=N\Delta s_{k},~e_{\min}\le e_N\le e_{\max}
\end{aligned}\right.
\end{aligned}
\end{equation}
\section{Interior Point Methods for Large Scale NLP} \label{sec:optim}
Although both sequential quadratic programming (SQP) and interior-point (IP) methods solve large-scale non-linear programming (NLP) problems \cite{wright1999numerical} like (\ref{eq:ocp}), recent benchmarks show the superior performance of IP solvers such as \texttt{IPOPT} \cite{wachter2006ipopt} and \texttt{KNITRO} \cite{byrd2006knitro} over SQP alternatives like \texttt{SNOPT} \cite{Mittelmann_AMPL_NLP_Benchmark_2025}, which lies in-line with our findings in Sec. \ref{sec:results}. However, the robust performance of IP methods is contingent upon several key factors: the strategy for handling non-convexity (i.e., an indefinite Hessian), the barrier parameter update rule, and the choice of globalization technique to ensure convergence \cite{wright1999numerical}. Therefore, we present a brief overview of Interior-Point (IP) theory to accompany our discussion in Sec. \ref{sec:results} on the various solver properties that help improve performance.

The optimization problem in (\ref{eq:ocp}) can be re-written as (\ref{eq:stacked_optim}), where ${\eta}\in\mathbb{R}^{N\times n_x+n_u\times(N-1)}$ is the new stacked decision variable and $c$ is the cost function. All equality constraints are stacked into $c_E$, representing general non-linear equality constraints. Similarly, inequality constraints (track/path bounds), (control, control rate bounds), (tire saturation) can be stacked into $c_I\in\mathbb{R}^{I_m}$ representing general nonlinear inequality constraints. Slack variables $s_i\in\mathbb{R}^{I_m}$, $s_i>0$ and a barrier parameter $\mu$ are introduced. An IP algorithm solves a set of approximate barrier problems for a sequence of positive barrier parameters $\mu_k$ that converge to zero \cite{byrd2006knitro}.
\begin{align}\label{eq:stacked_optim}
\min_{\eta,s}\quad c(\eta)\;-\;\mu \sum\nolimits_{i=1}^{I_m}\log s_i,\\
\text{s.t.}\quad c_E(\eta)=0,\quad
 c_I(\eta)-s=0\nonumber
\end{align}
we write the KKT conditions \cite{wright1999numerical} for problem \eqref{eq:stacked_optim} as:
\begin{align*}
    \nabla c(\eta) - A_E^T(\eta)y - A_I^T(\eta)z = 0, \quad & -\mu e + Sz = 0, \\
    c_E(\eta) = 0, \quad & c_I(\eta) - s = 0,
\end{align*}
where $y$ and $z$ are vectors of Lagrange multipliers, $e=(1, ..., 1)^T$, $S = \text{diag}(s_1, ..., s_m)$, $A_E$ and $A_I$ are the Jacobian matrices of constraints $c_E(\eta)$ and $c_I(\eta)$, respectively.  Defining the Lagrangian $\mathcal{L}$, the merit function for measuring progress towards a feasible and optimal solution $\phi$, and the corresponding Newton step to update the primal and dual variables $\eta, s, y, z$ as follows:
\begin{equation}
    \mathcal{L}(\eta, s, y, z) = c(\eta) - y^T c_E(\eta) - z^T(c_I(\eta) - s),
    \label{eq:Lagrangian}
\end{equation}
\begin{equation}
    \phi_{\nu}(\eta, s) = c(\eta) - \mu \sum_{i=1}^{m} \log s_i + \nu \|c_E(\eta)\|_2 + \nu \|c_I(\eta) - s\|_2,
    \label{eq:merit_func}
\end{equation} 
\begin{equation}
\begin{bmatrix}
    \nabla_{\eta\eta}^2 \mathcal{L} & 0 & -A_E^T(\eta) & -A_I^T(\eta) \\
    0 & Z & 0 & S \\
    A_E(\eta) & 0 & 0 & 0 \\
    A_I(\eta) & -I & 0 & 0
\end{bmatrix}
\begin{bmatrix}
    d_\eta \\ d_s \\ d_y \\ d_z
\end{bmatrix}
= 
\begin{bmatrix}
    A_E^T(\eta)y - \\\nabla f(\eta)+\\ A_I^T(\eta)z\\
    \mu e - Sz \\
    -c_E(\eta) \\
    s - c_I(\eta)
\end{bmatrix}\label{eq:newton_step}
\end{equation}
where $\nu > 0$. If the matrix in eq (\ref{eq:newton_step}) is well defined in terms of inertia, then the step $d$ is a descent direction for the merit function $\phi$. Through backtracking line search \cite{byrd2006knitro,wright1999numerical}, step-length $\alpha$ can be computed, and the decision variables at the next iterate $\eta^+, s^+, y^+, z^+$ can be obtained using:
\begin{align}
    \eta^{+} &= \eta + \alpha_s d_\eta, & s^{+} &= s + \alpha_s d_s,\\
    y^{+} &= y + \alpha_z d_y, & z^{+} &= z + \alpha_z d_z.
\end{align}
Our analysis in the results Sec. (\ref{sec:results}) demonstrate how the treatment of bad inertia or indefiniteness of the Hessian matrix $\nabla_{\eta\eta}^2\mathcal{L}$ in the Newton step (\ref{eq:newton_step}) and the barrier parameter ($\mu$) update strategy, dictate solver performance.
\section{Stiffness and Integrator Study}\label{sec:integrator_study}
Implicit and explicit integration techniques are known to perform well for stiff and non-stiff ODEs, respectively \cite{book_non_stiff, book_stiff}. To inform our choice of integrator, we first calculate the stiffness ratio at each node along a converged time-optimal trajectory for an oval track. Given the dynamics map $f$ in (\ref{eq:dynamics_full}), i.e. $\dot{x}(s) = f(x(s), u(s))$, the stiffness ratio $SR$ can be defined as the ratio of the magnitude of the fastest-decaying mode (largest real $\Re$ negative eigenvalue $\lambda_j$ magnitude) to the slowest-decaying mode (smallest real negative eigenvalue magnitude) of the Jacobian $J_k$ of $f$:
\begin{align}
\delta\dot x &= J_k\,\delta x,\quad J_k=\left.\frac{\partial f}{\partial x}\right|_{(x_k,u_k)},k=0,\dots,N\\
J_kv_j &= \lambda_j v_j,\ j=1,\ldots,n_x,\quad \mathcal S=\{\,j:\Re\lambda_j<0\,\},\\
r_j &= |\Re\lambda_j|,\ j\in\mathcal S,
\mathrm{SR}=\dfrac{\max_{j\in\mathcal S} r_j}{\min_{j\in\mathcal S} r_j},~if~
\mathcal \min_{j} r_j>0.\nonumber
\end{align}
We calculate a \textit{participation matrix} using the right and left eigenvalues for the Jacobians $J_k$ and backtrack the contribution of each state at every node in the stiffness ratio.

\begin{figure}[htbp]
  \centering
  \includegraphics[width=.92\columnwidth]{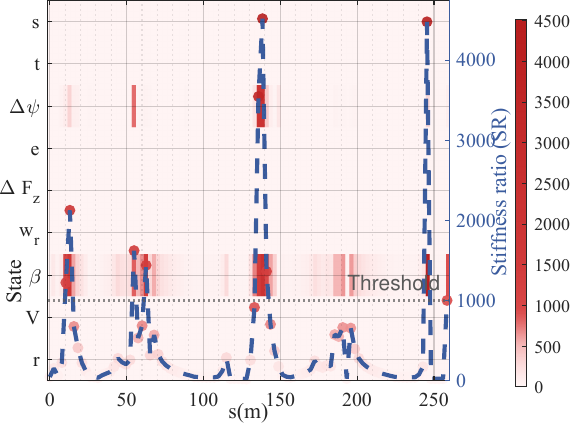}
  \caption{Node-wise stiffness ratio computed around a time-optimal race trajectory with backtracking state contributions.}\label{fig:stiff}
  \vspace{-3mm}
\end{figure}

\emph{The results in Figure \ref{fig:stiff} show the system is conditionally stiff, with the $SR$ peaking above $4500$ during aggressive maneuvers, far exceeding the traditional stiff threshold of $10^3$.} This stiffness originates primarily from the fast dynamics of the sideslip angle ($\beta$), which is consistent with the vehicle's physical behavior during high-speed cornering. This result implies the need to use implicit numerical integrators to prevent destabilization of iterations during optimization.

We study the impact of different integrators for the map $\mathcal{F}$ in (\ref{eq:ocp}) primarily on two metrics: (a) the solution accuracy compared to a ground truth, (b) the computational time required to solve the OCP (\ref{eq:ocp}). Each integrator is evaluated in combination with three solvers—\texttt{IPOPT}, \texttt{KNITRO}, and \texttt{SNOPT} \cite{gill2005snopt}. The test problem is a time-optimal $260$ m oval track trajectory ($N=100$), initialized from a dynamically feasible centerline guess. This study is motivated by the need to match an integrator's properties, such as its stability and order, with the underlying system dynamics and its characteristics, such as stiffness. We benchmark nine numerical integrators: explicit methods (RK2, RK4) and implicit methods (Implicit Euler, Backward Differentiation Formula (BDF) 4, 5, 6, Crank-Nicolson, Adams-Moulton 3-stage, and Gauss-Legendre 2-stage). A key property for differentiating implicit schemes can be A-stability; an A-stable integrator guarantees that the numerical solution to a stable physical problem will not become unstable, regardless of the integration step size \cite{book_non_stiff}. We note that among the methods tested, only the Implicit Euler, Crank-Nicolson, and Gauss-Legendre methods are A-stable. The ground truth solution was generated with a 6th-order A-stable Gauss-Legendre-3s method, verified against CVODES integration solver in Casadi \cite{casadi}. Although they produced identical results, both were omitted from the benchmark study due to their prohibitive solve times.

\section{Robust Initial Guess Formulation} \label{sec:initial_guess}
A strategy for initializing complex OCPs is first to generate a dynamically consistent trajectory to serve as an initial guess. While a dynamically feasible guess can help initialize the main OCP (\ref{eq:ocp}), generating this guess itself requires solving another OCP of similar complexity. This approach effectively layers one complex optimization problem on top of another, increasing the overall computational burden and introducing further points of failure. To circumvent the complexity of generating a full trajectory guess, we propose a simple yet powerful initialization method.

We propose a box-constrained nonlinear least square problem minimizing the dynamic rates (subset of the full system state) given by $\dot{x}_{\mathrm{dyn}} = [\dot{r}, \dot{V}, \dot{\beta}, \dot{\omega}_r, \dot{\Delta{F}}_z]^T$, which is solved to an optimality cost of zero, ensuring feasibility of the dynamic part of the equations of motion (\ref{eq:dynamics_full}).
We restrict ourselves to finding a single dynamically feasible set of state and control pairs for a given vehicle, making it computationally cheap, and initialize the dynamic components ($x_{\mathrm{dyn}}$) at every node of the OCP with this pair. Since direct multiple shooting is employed, that is every node has its own state and control variables, this computed state-control pair should satisfy the ${x_{\mathrm{dyn}}}_{k+1}=\mathcal{F}({x_{\mathrm{dyn}}}_{k+1},{x_{\mathrm{dyn}}}_k,u_k;\Delta s_k)$ at every node of the OCP  (as rates are zero) bringing down the feasibility error of the guess. The kinematic components $e,\Delta\phi$ are initialized with zeros representing centerline position, $s$ is known a priori, and $t$ is initialized by using an average velocity guess $V_{avg}$ and total path distance $s_{\mathrm{end}} - s_0$. The initial guess is formulated as follows:
\begin{equation}
\begin{aligned}
&\min_{\{x_{\mathrm{guess}}\},\{{u}_{\mathrm{guess}}\}} \dot{r}^2+\dot{V}^2+\dot{V}_{\omega_r}^2+\dot{\Delta F_z}^2 +\dot{\beta}^2 + (V-V_{\omega_r})^2 \\[-1pt]
& \text{\scriptsize(bounds)}\ \left\{
\begin{aligned}
& 0^-\le r_{\mathrm{guess}}\le 0^+, \\&V_{\min}\le V_{\mathrm{guess}}\le V_{\max}, \\&\underline{u}\ \le\ {u}_{\mathrm{guess}}\ \le\ \overline{u}
\end{aligned}\right.\\
\end{aligned}\label{eq:guess_ocp}
\end{equation}
We emphasize that the cost function and bounds in the above problem can be customized to produce different responses; we opted to obtain a straight line driving condition and introduce $(V-V_{\omega_r})^2$ to minimize slip. The guess once computed for a vehicle can be fixed and \textit{is not required to be recomputed}. 

\section{Results}\label{sec:results}
We frame the OCP (\ref{eq:ocp}) for a Lexus LC 500 vehicle model\cite{djeumou2025reference} in MATLAB R2024b via CasADi \cite{casadi} interfaced with \texttt{IPOPT}\cite{wachter2006ipopt}, \texttt{SNOPT}\cite{gill2005snopt}, and \texttt{KNITRO}\cite{byrd2006knitro} NLP solvers. The benchmarks were run on a desktop PC with a 5.7GHz AMD Ryzen 9 9950X processor, with Just-In-Time (JIT) compilation enabled for all results except those in Table \ref{tab:solver_times}.

\subsection{Integrator and Solver Tandem Study Results}
\begin{figure*}[t]
  \centering
  \includegraphics[width=1\textwidth]{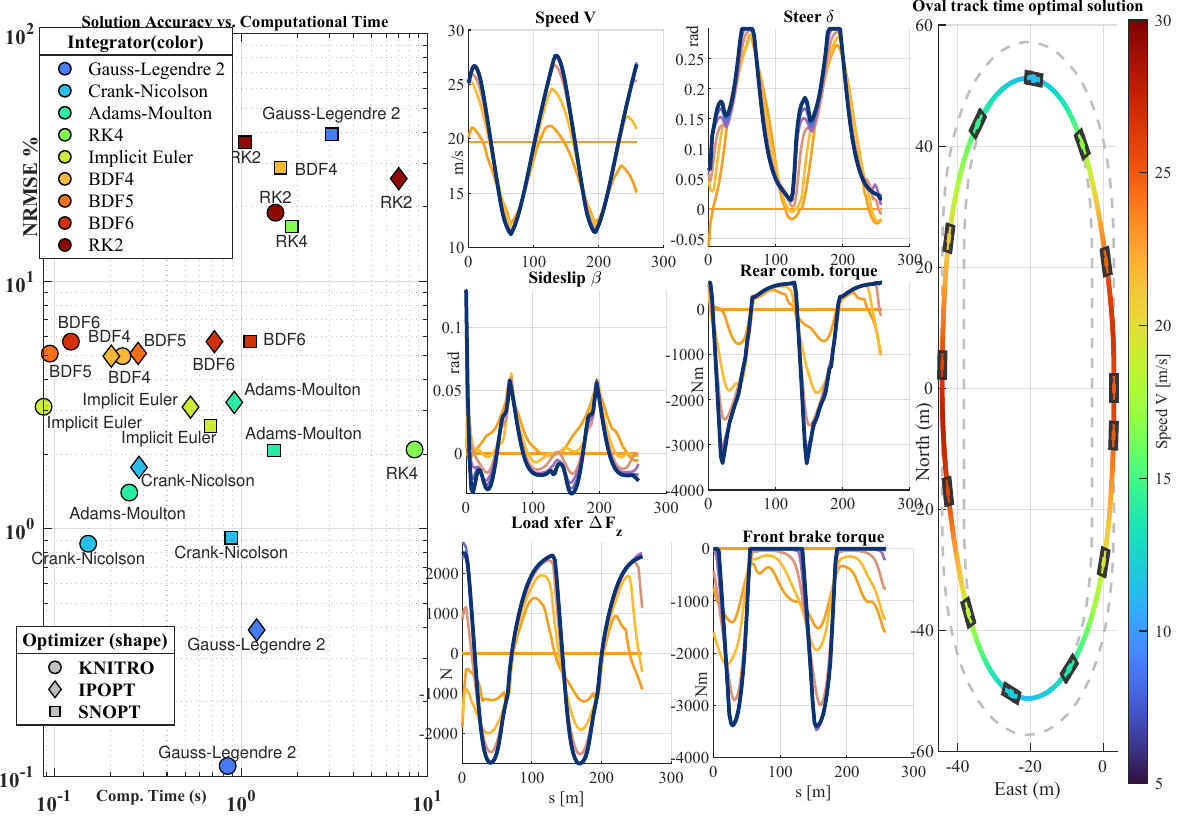}
  \caption{(a) Different numerical integrators and optimization solver pairs plotted with NRMSE$\%$ depicting accuracy (Y-axis) and computation time (X-axis). (b) Iteration-wise converged time-optimal solution generated under $55ms$ $(22)$ iterations for an oval track starting from the introduced guess strategy, illustrating optimal vehicle path, states, and control sequences.}
  \label{fig:integ_n_oval}
  \vspace{-5mm}
\end{figure*}

\begin{table}[htbp]
\centering
\caption{CPU times (s) by integrator and solver. \textbf{Bold} = fastest optimal; {\color{red}red} = suboptimal ; {\color{red}NS} = No Solution.}
\label{tab:solver_times}
\begin{tabular}{lccc}
\toprule
Integrator           & KNITRO  & IPOPT   & SNOPT   \\
\midrule
Implicit Euler       & \textbf{0.087}       & 0.535   & 0.686   \\
BDF5                 & \textbf{0.094}       & 0.280   & {\color{red}NS}       \\
BDF6                 & \textbf{0.122}       & 0.719   & 1.120   \\
Crank–Nicolson      & \textbf{0.151}       & 0.283   & 0.885   \\
BDF4                 & 0.231       & \textbf{0.201}       & {\color{red}1.630}   \\
Adams–Moulton 3s       & \textbf{0.251}       & 0.919   & 1.500   \\
Gauss–Legendre 2s     & \textbf{0.847}       & 1.212   & {\color{red}3.070}   \\
RK4                  & \textbf{8.539}       & {\color{red}NS}       & {\color{red}1.870}   \\
RK2                  & {\color{red}1.530}   & {\color{red}7.024}   & {\color{red}1.050}   \\
\bottomrule
\end{tabular}
\end{table}
The comprehensive results for the integrator study introduced in Sec. \ref{sec:integrator_study} are summarized in Figure \ref{fig:integ_n_oval}, and Table \ref{tab:solver_times}. Solution accuracy is quantified by the normalized root mean square error (NRMSE) for each control input $u_j,j=1,2,3$, which is the standard root mean square error between the computed trajectory $u_j$ and the ground truth, $u_j^{*}$. Figure \ref{fig:integ_n_oval} illustrates the trade-off between solution accuracy and computational cost by plotting the maximum NRMSE \% made in a control input against the CPU solve time for each integrator-solver pair (averaged over 10 consecutive solves of each combination). Table \ref{tab:solver_times} provides a detailed breakdown of the CPU times for each integrator-solver pair. Note that rather than employing hidden internal Newton iterations to resolve the implicit dynamics at each step, the optimizer solves for the states, controls, and any multi-stage variables simultaneously as a single NLP problem. We summarize our findings as follows:\\
\noindent\textbullet~ \textbf{Implicit vs Explicit methods:}
The results confirm the unsuitability of explicit integrators for this problem, aligning directly with our stiffness analysis (Sec. \ref{sec:integrator_study}). The stiff dynamics impose severe penalties on explicit schemes, placing both RK2 and RK4 in the top-right quadrant of Figure \ref{fig:integ_n_oval}(a), characterized by high computational cost and poor accuracy. None of the solvers converge to the optimal solution when using the single-stage RK2 method, as shown in Table \ref{tab:solver_times}. Although KNITRO succeeds with the four-stage RK4, its solve time of $8.54 $ seconds is prohibitive for real-time use as it is nearly \textit{$100$ times slower} than the implicit Euler scheme.

\noindent\textbullet~ \textbf{Accuracy and A-stability:}  A critical insight from Figure \ref{fig:integ_n_oval} (a) is that for this stiff dynamics, simply using a higher order integrator does not guarantee a more accurate solution. Instead, the property of A-stability emerges as the dominant factor. A clear pattern is visible among the A-stable implicit methods. As the order increases from the 1st-order implicit Euler ($3$--$5\%$  error) to the 2nd-order Crank-Nicholson ($<1\%$  error), and finally to the 4th-order Gauss-Legendre 2 ($0.1\%$  error), the solution accuracy consistently improves. Conversely, non-A-stable, higher-order methods like BDF 4-6 and Adams-Moulton 3-stage violate this trend, yielding larger errors than the second-order Crank-Nicholson.

\noindent\textbullet~ \textbf{Computational Efficiency:} To identify the best configuration for our problem, we fix an upper bound on the maximum NRMSE $\%$ to be $5\%$, and the upper bound on computational time to be $0.2$ seconds. Within these bounds, the combination of the A-stable implicit Euler integrator and the KNITRO solver performs best, achieving convergence in just $0.087$ seconds (Table \ref{tab:solver_times}) with an error margin of $3-4\%$.

\noindent\textbullet~ \textbf{Solver Performance:} As shown in Table \ref{tab:solver_times} and Figure \ref{fig:integ_n_oval}, the second-order IP solvers (KNITRO, IPOPT) exhibit superior robustness than the first-order SQP-based SNOPT, which frequently converges to suboptimal solutions. This performance gap underscores the importance of exact Hessian information for reliable convergence in highly nonlinear problems, where the quasi-Newton approximations employed by SNOPT appear insufficient. While both IP solvers demonstrate consistent robustness, KNITRO achieves notably faster solve times—often two to three times faster than IPOPT.


\subsection{Robust Initial Guess Results}
\begin{table}[htbp]
\caption{KNITRO per-segment performance: zero guess vs. robust initial guess. \;{\small {\color{red}NS} = evaluation/convergence failure}}
\centering
\footnotesize
\setlength{\tabcolsep}{3pt}
\renewcommand{\arraystretch}{1.05}
\begin{tabular}{@{}lcc|cc|c@{}}
\toprule
 & \multicolumn{2}{c|}{Zero guess} & \multicolumn{2}{c|}{Robust init} & \multirow{2}{*}{Comp time gain $(\times)$} \\
\cmidrule(lr){2-3}\cmidrule(lr){4-5}
Segment & Iter & Time [s] & Iter & Time [s] &  \\
\midrule
S01 & \textcolor{red}{NS} & \textcolor{red}{NS} & 27 & 0.073 & \textcolor{red}{NS} \\
S02 & 2005 & 7.674 & 40 & 0.098 &  $\times$78.30 \\
S03 & \textcolor{red}{NS} & \textcolor{red}{NS} & 40 & 0.105 & \textcolor{red}{NS} \\
S04 & 807  & 2.967 & 62 & 0.153 &  $\times$19.39 \\
S05 & \textcolor{red}{NS} & \textcolor{red}{NS} & 62 & 0.139 & \textcolor{red}{NS} \\
S06 & 1227 & 5.007 & 35 & 0.085 &  $\times$58.90 \\
S07 & 693  & 2.585 & 39 & 0.100 &  $\times$25.85 \\
S08 & 202  & 0.689 & 66 & 0.159 &  $\times$4.34 \\
S09 & 1675 & 6.542 & 38 & 0.091 &  $\times$71.89 \\
S10 & 937  & 3.763 & 49 & 0.136 &  $\times$27.67 \\
\bottomrule
\end{tabular}
\label{tab:guess_study}
\end{table}
Using the optimal integrator-solver pair (Implicit Euler with \texttt{KNITRO}) we found for our problem, we evaluate the impact of our robust initialization strategy (Sec. \ref{sec:initial_guess}) against an ill-informed, zero guess. The feasibility error representing the maximum constraint violation is computed for the time-optimal OCP (\ref{eq:ocp}) around an oval track at the two guesses. The robust guess provides a starting point with an error of just $9\times 10^{-1}$; conversely, the ill-informed guess is strongly infeasible, beginning with an error of $6.7\times 10^3$. Moreover, when the same time-optimal problem is solved across the following robustness test for $10$ different racing track segments \cite{thompson2024adaptive}, a similar difference in initial feasibility error is observed. This \textit{four-order-of-magnitude decrease} in initial infeasibility effectively transforms the optimization landscape, providing a better starting point for the optimizer. Figure \ref{fig:integ_n_oval} (b) illustrates the iteration-wise descent (light yellow to dark blue) towards the optimal solution that took only $55$ ms and $22$ iterations to converge. Particularly, the dynamic states and controls are initialized with the same constant $x_{\mathrm{guess}},u_{\mathrm{guess}}$ computed using the proposed guess formulation in Sec. \ref{sec:initial_guess}, from where the optimization proceeds towards an optimum.

Further, to isolate the impact of the robust initial guess from the robustness of the optimizer, keeping the integrator-solver pair intact, we solve a time-optimal OCP on $10$ distinct $250$ m ($N=100$) segments of a racetrack \cite{thompson2024adaptive}, starting from the fixed robust guess and a zero guess (initializing all decision variables with $0$). The results, presented in Table \ref{tab:guess_study}, highlight three critical findings. First, the robust guess is essential for convergence. The zero guess fails to find a solution in 30\% of the cases (S01, S03, S05). Second, the robust guess enables real-time performance. The zero guess never meets the $200$ ms real-time target, with its fastest solve time being around $689$ ms, whereas the robust guess consistently delivers solutions in under $160$ ms. Finally, the proposed strategy is both efficient and generalizable; the same guess once computed works on all $10$ diverse segments.

\subsection{Optimization Solvers}
\begin{figure}[!htbp]
  \centering
  \includegraphics[width=0.9\columnwidth]{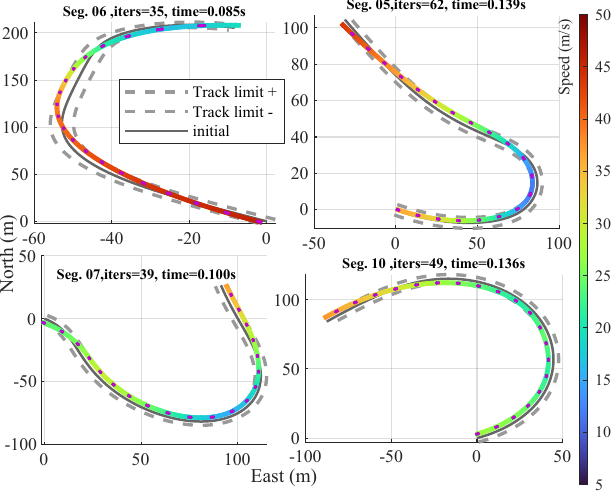}
  \caption{Racing track segments time-optimal paths.}
  \label{fig:thill_plots}
  \vspace{-3mm}
\end{figure}
\begin{figure}[htbp]
  \centering
  \includegraphics[width=.8\columnwidth]{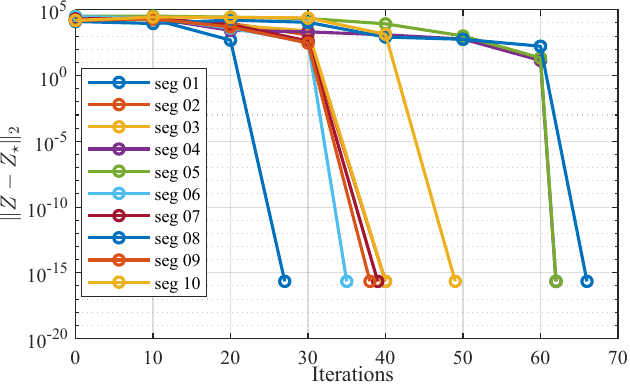}
  \caption{Racing track segments convergence sequence.}
  \label{fig:conv_sequence}
\end{figure}

Since \texttt{SNOPT} consistently produced sub-optimal solutions (Table \ref{tab:solver_times}) in the integrator test study, we limit our discussion to the parameters that are imperative to stabilize and produce better performance for the IP Methods.\\
\noindent\textbullet~  \textbf{Barrier Parameter Influence:} The barrier parameter update strategy is a critical hyperparameter in IP methods, directly impacting convergence and solution quality \cite{byrd2006knitro,nocedal2009adaptive}. Our tests show both IP solvers are highly sensitive to this choice, with a poor selection leading to suboptimal or failed solutions (Table \ref{tab:barmurule-summary}). For \texttt{IPOPT}, we found its default monotone strategy frequently failed to converge across different tests, in contrast to the adaptive strategy, which was found to be more robust. A similar effect was observed for \texttt{KNITRO}, where the advanced dampmpc strategy—a safeguarded predictor-corrector rule—solved the problem in just $22$ iterations, while a basic monotone approach failed (Table \ref{tab:barmurule-summary}). Further details on various barrier strategies can be found in \cite{nocedal2009adaptive} (excluded here for brevity).
\begin{table}[htbp]
\centering
\caption{\texttt{KNITRO} barrier update strategies with iteration count for oval OCP}
\footnotesize 
\setlength{\tabcolsep}{2pt}
\renewcommand{\arraystretch}{1.05}
\begin{tabular}{@{}lcccccc@{}}
\toprule
 & dampmpc & probing & fullmpc & adaptive & quality & monotone \\
\midrule
Iterations & \textbf{22} & 29 & 29 & 35 & 50 & 289 \\
Optimality & optimal & optimal & optimal & optimal & optimal & \textcolor{red}{Sub-optimal} \\
\bottomrule
\end{tabular}
\label{tab:barmurule-summary}
\end{table}\\
\noindent\textbullet~  \textbf{Solver Robustness, Real-Time Performance, and Convergence:} Building on our prior results, we next compare \texttt{IPOPT} and \texttt{KNITRO} on $10$ time-optimal distinct race track segments. We use our robust guess to evaluate each solver's capacity for real-time convergence ($<200$ ms), reporting the results in Table \ref{tab:robust-ipopt-vs-knitro}. The benchmark results underscore \texttt{KNITRO}'s robustness for real-time trajectory optimization, as it consistently meets the $200$ ms performance target across all segments. In contrast, \texttt{IPOPT} proves unreliable, exceeding the time limit on $90\%$ of the cases. This performance gap is starkly highlighted by segment $S06$, where \texttt{IPOPT}'s solve time exceeds $10$ seconds (due to a sharp increase in per-iteration cost)—over $126$ times slower than \texttt{KNITRO}—rendering it unsuitable for time-critical tasks. Across all ten segments in Figure \ref{fig:conv_sequence}, \texttt{KNITRO} rapidly drives the optimality error \( \lVert Z - Z_\star \rVert_2 \) down by multiple orders of magnitude and then settles into a uniform linear-rate tail, reaching tolerance by \( < 70 \) iterations without oscillations. A few of the converged time-optimal path trajectories obtained from \texttt{KNITRO} in Table \ref{tab:robust-ipopt-vs-knitro} are visualized in Figure \ref{fig:thill_plots}.


\begin{table}[htbp]
\caption{Robust guess with IPOPT (adaptive $\mu$) $\&$ KNITRO (dampmpc $\mu$). Times $> \,$\textbf{0.2 s} in \textcolor{red}{red}, least in \textbf{bold}}
\centering
\footnotesize
\setlength{\tabcolsep}{3pt}
\renewcommand{\arraystretch}{1.05}
\begin{tabular}{@{}lccccc@{}}
\toprule
 & \multicolumn{2}{c}{IPOPT (adaptive $\mu$)} & \multicolumn{2}{c}{KNITRO (dampmpc $\mu$)} & \multirow{2}{*}{Time gain $(\times)$} \\
\cmidrule(lr){2-3}\cmidrule(lr){4-5}
Segment & Time [s] & Iter & Time [s] & Iter & \\
\midrule
S01 & \textcolor{red}{0.406} & 82  & \textbf{0.073} & 27 & $\times$5.56 \\
S02 & 0.178                   & 36  & \textbf{0.098} & 40 & $\times$1.82 \\
S03 & \textcolor{red}{0.415} & 102 & \textbf{0.105} & 40 & $\times$3.95 \\
S04 & \textcolor{red}{0.485} & 111 & \textbf{0.153} & 62 & $\times$3.17 \\
S05 & \textcolor{red}{0.461} & 101 & \textbf{0.139} & 62 & $\times$3.32 \\
S06 & \textcolor{red}{10.714} & 129 & \textbf{0.085} & 35 & $\times$126.05 \\
S07 & \textcolor{red}{0.783} & 90  & \textbf{0.100} & 39 & $\times$7.83 \\
S08 & \textcolor{red}{0.323} & 76  & \textbf{0.159} & 66 & $\times$2.03 \\
S09 & \textcolor{red}{0.324} & 61  & \textbf{0.091} & 38 & $\times$3.56 \\
S10 & \textcolor{red}{0.570} & 81  & \textbf{0.136} & 49 & $\times$4.19 \\
\bottomrule
\end{tabular}
\label{tab:robust-ipopt-vs-knitro}
\end{table}

\noindent\textbullet~  \textbf{\texttt{KNITRO} and \texttt{IPOPT}: Globalization and Bad Inertia Strategies:} We specifically pay attention to the case of $S06$ where \texttt{IPOPT} takes about two orders of magnitude more time as compared to \texttt{KNITRO}, for just $129$ iterations, attributing to a transient blow-up in per-iteration linear-algebra cost. Both solvers employ different globalization techniques, such as using a merit function like (\ref{eq:merit_func}) in \texttt{KNITRO} or a filter-based method \cite{wachter2006ipopt}, that can influence how the iterates proceed. However, the primary difference stems from handling bad inertia in the Newton step or indefiniteness of the Hessian $\nabla_{\eta\eta}^2 \mathcal{L}$ eq (\ref{eq:newton_step}). On segment $S06$, \texttt{IPOPT}’s time surge stems from its inertia-correction strategy: when the KKT system is indefinite or nearly singular, IPOPT repeatedly re-factorizes with regularization to enforce inertia. In contrast, \texttt{KNITRO}’s hybrid IP algorithm safeguards by switching to a trust-region based conjugate gradient (CG) step when bad inertia is detected, stabilizing the step in a few CG iterations and achieving real-time performance.
\subsection{Emergency Collision Avoidance Test}
\begin{figure}[!t]
  \centering
  \includegraphics[width=.8\columnwidth]{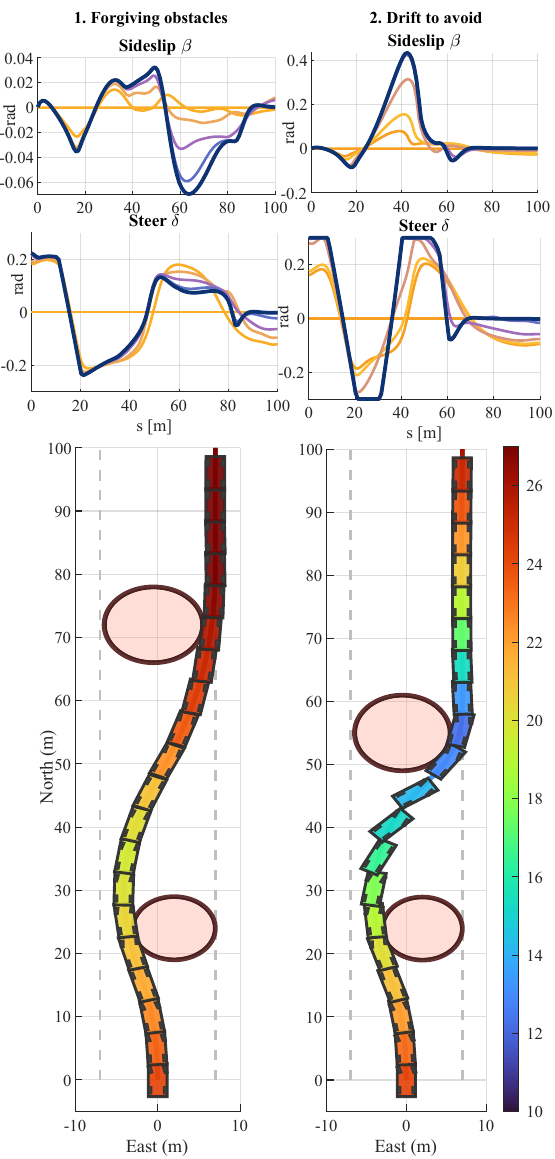}
  \caption{Collision avoidance test trajectories: left (forgiving obstacle), right (requires drifting maneuver).}
  \label{fig:col_avoid_path}
\end{figure}

We put our architecture to test on two high-speed ($90$~km/h) emergency collision avoidance scenarios similar to those in \cite{CA_scenarios}. Figure \ref{fig:col_avoid_path} highlights the framework's versatility, showing it can autonomously generate either a stable, low-sideslip maneuver ($<58$ ms) for a forgiving scenario or a controlled drift ($>0.4$~rad sideslip, $<72$~ms) for a more critical one. We empirically prove this aggressive drift is the \textit{only} viable solution, as reinstating the tire saturation constraints (\ref{eq:ocp}) renders the problem infeasible. Despite the complexity, both trajectories were generated in under $72$ ms, confirming real-time performance. This showcases our framework's ability to act as a unified planner that obviates the need for specialized drift controllers while providing a real-time \textit{certificate of feasibility}. 

Further, we generated similar high-speed oval racing and collision avoidance trajectories for a test all-terrain vehicle (ATV) to validate real-world plausibility. A full-fidelity model of this vehicle is simulated by BeamNG, and the optimal trajectories are generated using a low-fidelity, data-driven model (similar to (\ref{eq:dynamics_full})) trained on BeamNG data. We track the reference trajectory using a short-horizon Model Predictive Control (MPC) running online.

\begin{figure}[tb]
  \centering
  \includegraphics[width=.8\columnwidth]{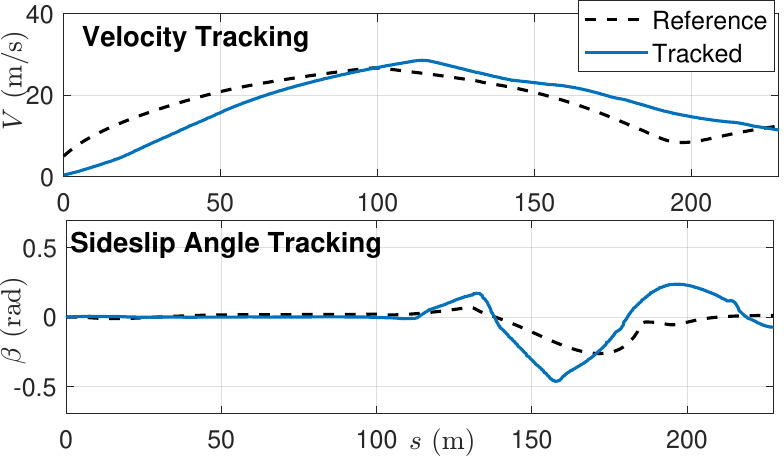}
  \caption{Tracking performance velocity and sideslip for high-speed emergency collision avoidance in BeamNG.}
  \label{fig:beam_states}
\end{figure}

Figure \ref{fig:beamNG} and Figure \ref{fig:beam_states} show successful double obstacle avoidance at high velocities $90$ km/hr. The performance plots (Figure \ref{fig:beam_states}) show an initial velocity error due to a standing start and a maximum sideslip tracking error of up to $\pm15$ deg. The sideslip deviation is a deliberate consequence of the online tracking MPC's tuning, which prioritized path tracking over minimizing transient dynamic state errors to ensure safety. Future work aims to test these trajectories on the physical model of this ATV vehicle under development.

\section{CONCLUSIONS}

We demonstrate that real-time, reference-free trajectory generation at the limits of vehicle handling is achievable through a principled initial guess strategy, integrator selection informed by system stiffness analysis, and tailored interior point solver techniques for barrier updates and nonconvexity handling. The framework was validated across multiple race track segments for time-optimal planning, and in a high-speed collision avoidance scenario that requires drift for feasibility. Future work will implement the proposed architecture on real-world race cars. We will investigate how the current stiffness analysis and equilibrium-based initial guess generalizes to other system dynamics.

\bibliographystyle{IEEEtran}
\bibliography{bib}

@article{CA_scenarios,
title = {Planning and control of drifting-based collision avoidance strategy under emergency driving conditions},
journal = {Control Engineering Practice},
volume = {139},
pages = {105625},
year = {2023},
issn = {0967-0661},
author = {Daofei Li and Jiajie Zhang and Siyuan Lin}
}

@misc{beamng_tech,
    title = "{BeamNG.tech}",
    author = {{BeamNG GmbH}},
    howpublished = {\url{https://www.beamng.tech/}},
    year = {2025},
    note = {Version 0.35.0.0, Accessed on June 15, 2025}
}

@article{gill2005snopt,
  title={SNOPT: An SQP algorithm for large-scale constrained optimization},
  author={Gill, Philip E and Murray, Walter and Saunders, Michael A},
  journal={SIAM review},
  volume={47},
  number={1},
  pages={99--131},
  year={2005},
  publisher={SIAM}
}

@misc{Mittelmann_AMPL_NLP_Benchmark_2025,
  author       = {Hans Mittelmann},
  title        = {AMPL{-}NLP Benchmark},
  year         = {2025},
  month        = jun,
  institution  = {Arizona State University},
  note         = {Updated June 28, 2025. Benchmarks for Optimization Software},
  urldate      = {2025-09-13}
}

@Article{casadi,
  author = {Joel A E Andersson and Joris Gillis and Greg Horn
            and James B Rawlings and Moritz Diehl},
  title = {{CasADi} -- {A} software framework for nonlinear optimization
           and optimal control},
  journal = {Mathematical Programming Computation},
  volume = {11},
  number = {1},
  pages = {1--36},
  year = {2019},
  publisher = {Springer}
}

@article{wachter2006ipopt,
  title={On the implementation of an interior-point filter line-search algorithm for large-scale nonlinear programming},
  author={W{\"a}chter, Andreas and Biegler, Lorenz T},
  journal={Mathematical programming},
  volume={106},
  number={1},
  pages={25--57},
  year={2006},
  publisher={Springer}
}

@article{wright1999numerical,
  title={Numerical optimization},
  author={Wright, Stephen and Nocedal, Jorge and others},
  journal={Springer Science},
  volume={35},
  number={67-68},
  pages={7},
  year={1999}
}

@incollection{byrd2006knitro,
  title={KNITRO: An integrated package for nonlinear optimization},
  author={Byrd, Richard H and Nocedal, Jorge and Waltz, Richard A},
  booktitle={Large-scale nonlinear optimization},
  pages={35--59},
  year={2006},
  publisher={Springer}
}

@article{nocedal2009adaptive,
  title={Adaptive barrier update strategies for nonlinear interior methods},
  author={Nocedal, Jorge and W{\"a}chter, Andreas and Waltz, Richard A},
  journal={SIAM Journal on Optimization},
  volume={19},
  number={4},
  pages={1674--1693},
  year={2009},
  publisher={SIAM}
}

@book{talbot2024optimal,
  title={Optimal vehicle control under friction uncertainty-from driver assistance to drift transitions},
  author={Talbot, John Andrew},
  year={2024},
  publisher={Stanford University}
}

@inproceedings{domberg2022deep,
  title={Deep drifting: Autonomous drifting of arbitrary trajectories using deep reinforcement learning},
  author={Domberg, Fabian and Wembers, Carlos Castelar and Patel, Hiren and Schildbach, Georg},
  booktitle={2022 International Conference on Robotics and Automation (ICRA)},
  pages={7753--7759},
  year={2022},
  organization={IEEE}
}

@article{goh2024beyond,
  title={Beyond the stable handling limits: nonlinear model predictive control for highly transient autonomous drifting},
  author={Goh, Jonathan YM and Thompson, Michael and Dallas, James and Balachandran, Avinash},
  journal={Vehicle System Dynamics},
  volume={62},
  number={10},
  pages={2590--2613},
  year={2024},
  publisher={Taylor \& Francis}
}

@article{weber2023modeling,
  title={Modeling and control for dynamic drifting trajectories},
  author={Weber, Trey P and Gerdes, J Christian},
  journal={IEEE Transactions on Intelligent Vehicles},
  volume={9},
  number={2},
  pages={3731--3741},
  year={2023},
  publisher={IEEE}
}

@inproceedings{djeumou2025reference,
  title={Reference-free formula drift with reinforcement learning: From driving data to tire energy-inspired, real-world policies},
  author={Djeumou, Franck and Thompson, Michael and Suminaka, Makoto and Subosits, John},
  booktitle={2025 IEEE International Conference on Robotics and Automation (ICRA)},
  pages={3610--3616},
  year={2025},
  organization={IEEE}
}

@book{goh2019_thesis,
  title={Automated vehicle control beyond the stability limits},
  author={Goh, Jonathan Yan Ming},
  year={2019},
  publisher={Stanford University}
}

@article{thompson2024adaptive,
  title={Adaptive nonlinear model predictive control: maximizing tire force and obstacle avoidance in autonomous vehicles},
  author={Thompson, Michael and Dallas, James and Goh, Jonathan YM and Balachandran, Avinash},
  journal={IEEE Transactions on Field Robotics},
  year={2024},
  publisher={IEEE}
}

@article{subosits2019racetrack,
  title={From the racetrack to the road: Real-time trajectory replanning for autonomous driving},
  author={Subosits, John K and Gerdes, J Christian},
  journal={IEEE Transactions on Intelligent Vehicles},
  volume={4},
  number={2},
  pages={309--320},
  year={2019},
  publisher={IEEE}
}

@article{chris_g_cascaded_22,
  title={Long-horizon vehicle motion planning and control through serially cascaded model complexity},
  author={Laurense, Vincent A and Gerdes, J Christian},
  journal={IEEE Transactions on Control Systems Technology},
  volume={30},
  number={1},
  pages={166--179},
  year={2021},
  publisher={IEEE}
}

@inproceedings{long_horizon_point_mass,
  title={Predictive control of autonomous ground vehicles with obstacle avoidance on slippery roads},
  author={Gao, Yiqi and Lin, Theresa and Borrelli, Francesco and Tseng, Eric and Hrovat, Davor},
  booktitle={Dynamic systems and control conference},
  volume={44175},
  pages={265--272},
  year={2010}
}

@inproceedings{motion_primitives,
  title={Predictive control for agile semi-autonomous ground vehicles using motion primitives},
  author={Gray, Andrew and Gao, Yiqi and Lin, Theresa and Hedrick, J Karl and Tseng, H Eric and Borrelli, Francesco},
  booktitle={2012 American Control Conference (ACC)},
  pages={4239--4244},
  year={2012},
  organization={IEEE}
}

@inproceedings{switch_logic,
  title={Collision avoidance with transitional drift control},
  author={Zhao, Tong and Yurtsever, Ekim and Chladny, Ryan and Rizzoni, Giorgio},
  booktitle={2021 IEEE International Intelligent Transportation Systems Conference (ITSC)},
  pages={907--914},
  year={2021},
  organization={IEEE}
}

@article{justifying_emergency,
  title={Justifying emergency drift control for automated vehicles},
  author={Zhao, Tong and Yurtsever, Ekim and Rizzoni, Giorgio},
  journal={IFAC-PapersOnLine},
  volume={55},
  number={24},
  pages={141--148},
  year={2022},
  publisher={Elsevier}
}

@book{book_stiff,
  title={Solving Ordinary Differential Equations II: Stiff and Differential-Algebraic Problems},
  author={Hairer, E. and N{\o}rsett, S.P. and Wanner, G.},
  isbn={9783540604525},
  lccn={86031456},
  series={Solving Ordinary Differential Equations II: Stiff and Differential-algebraic Problems},
  year={1993},
  publisher={Springer}
}

@book{book_non_stiff,
  title={Solving Ordinary Differential Equations II: Stiff and Differential-Algebraic Problems},
  author={Hairer, Ernst and Wanner, Gerhard},
  year={1996},
  publisher={Springer Berlin Heidelberg},
  series={Springer Series in Computational Mathematics},
  isbn={978-3-642-05220-0}
}

@article{Lund,
  title={An investigation of optimal vehicle maneuvers for different road conditions},
  author={Olofsson, Bj{\"o}rn and Lundahl, Kristoffer and Berntorp, Karl and Nielsen, Lars},
  journal={IFAC Proceedings Volumes},
  volume={46},
  number={21},
  pages={66--71},
  year={2013},
  publisher={Elsevier}
}

\end{document}